\documentstyle[a4,12pt,twoside]{article}

 \title{{\Large\bf Koszul Property and Poincar\'e Series \\ of Matrix Bialgebras of Type $A_n$ } \footnote{  JOURNAL OF ALGEBRA Volume 192,
N2, (1997), P. 734-748.}}
\author{\large  Phung Ho Hai\thanks{Suppoted by DFG at the Graduiertenkolleg ``Mathematik im Bereich ihrer Wechselwirkung mit der Physik'', Mathematisches Institut, University of Munich}\\
\small Sektion Physik der Universit\"at M\"unchen, Lehrstuhl Wess\\
\small Theresienstr 37,  80333 Munich, FRG\\ 
\small e-mail  phung@lswes8.ls-wess.physik.uni-muenchen.de}

\date{}
\begin{document}
\maketitle
\def\iii{^{-1}}
\def\ii{_{-1}}
\def\AAa{$A_\al ,\al\in J $}
\def\lora{\longrightarrow}
\def\ot{\otimes}
\def\opl{\oplus}
\def\loma{\longmapsto}
\def\Vn{{V^{\ot n}}}
\def\Vm{{V^{\ot m}}}
\def\Vsn{{V^{*\ot n}}}
\def\bigopinf{\bigoplus_{n=0}^\infty}
\def\al{\alpha}
\def\si{\sigma}
\newcommand{\bbas}{\begin{eqnarray*}}
\newcommand{\eeas}{\end{eqnarray*}}
\newcommand{\bbs}{\begin{displaymath}}
\newcommand{\ees}{\end{displaymath}}
\def\bb{\begin{equation}}
\def\ee{\end{equation}}
\def\eea{\end{eqnarray}}
\def\ochi{\overline{\chi}}
\def\taubar{\overline{\tau}}
\def\wbar{\overline{w}}
\def\KK{\mbox{$\mathbf K$}}
\def\ZZ{\mbox{$\mathbf Z$}}

\def\vn{{v_n}}
\def\vi{{v_i}} 
\def\vk{{v_k}}
\def\wk{{w_k}}
\def\Tvi{T_\vi}
\def\Tvk{T_\vk}
\def\Twk{T_\wk}
\def\En{E_1^{\ot n}}
\def\Esn{E_1^{*\ot n}}
\newcommand{\dayso}[1]{\mbox{$1,2,\cdots ,#1$}}
\renewcommand{\day}[2]{\mbox{$#1_1,#1_2,\cdots ,#1_{#2}$}}
\newcommand{\dayup}[2]{\mbox{$#1^1,#1^2,\cdots ,#1^{#2}$}}
\newcommand{\daysub}[2]{\mbox{$#1_1,#1_2,\cdots ,#1_{#2}$}}
\def\bba{\begin{eqnarray}}
\newtheorem{thm}{Theorem}[section]
\newtheorem{lem}[thm]{Lemma}
\newtheorem{dn}[thm]{Definition}
\newtheorem{bem}[thm]{Bemerkung}
\newtheorem{edl}[thm]{Theorem}
\newtheorem{dl}[thm]{Satz}
\newtheorem{cor}[thm]{Corollar}
\newtheorem{kor}[thm]{Korollar}
\newtheorem{pro}[thm]{Proposition}
\newtheorem{dnsatz}[thm]{Definition-Satz}
\def\Lambd{{\mit\Lambda}}
\def\ybo{Yang-Baxter-Operator}
\def\haa{Hecke-Algebra}
\def\yb{Yang-Baxter}
\def\ho{Hecke-Operator}
\def\ccc{,\cdots ,}
\def\RR{\mathbf{R}}
\def\CC{\mathbf{C}}
\def\det{\mbox{\rm det}}
\def\detq{\mbox{\rm det$_q$}}
\def\sdet{\mbox{\rm Ber$_q$}}
\newcommand{\Ai}{\mbox{$A^{-1}$}}
\newcommand{\Vi}{\mbox{$V^{-1}$}}
\newcommand{\Hi}{\mbox{$H^{-1}$}}
\newcommand{\Di}{\mbox{$D^{-1}$}}
\def\hi{\hat{i}}
\def\hj{\hat{j}}
\def\hm{\hat{m}}
\def\hmm{\hat{1}}
\def\hn{\hat{n}}
\def\hk{\hat{k}}
\def\hl{\hat{l}}
\def\hp{\hat{p}}
\def\hq{\hat{q}}
\def\ha{\hat{a}}
\def\hb{\hat{b}}
\def\hhh{\hat{2}}
\def\hs{\hat{s}}
\def\Ker{\mbox{\rm Ker}}
\def\Im{\mbox{\rm Im}}
\def\coend{\mbox{\rm coend}}
\def\Mor{\mbox{\rm Mor}}
\def\Hom{\mbox{\rm Hom}}
\def\Vek{\mbox{\rm Vek}}
\def\vek{\mbox{\rm vek}}
\def\vec{\mbox{\rm vec}}
\def\Vec{\mbox{\rm Vec}}
\def\Ber{\mbox{\rm Ber}}
\def\End{\mbox{\rm End}}
\def\kk{$k$}
\def\d{\mbox{\rm d}}
\def\kkk{$k$-lineare Basis}
\def\H{{\cal H}}
\def\W{{\cal W}}
\def\A{{\cal A}}
\def\B{{\cal B}}
\def\L{{\cal L}}
\def\C{\mbox{$\cal C$}}
\def\Z{{\cal Z}}
\def\S{{\cal S}}
\def\V{{\cal V}}
\def\P{{\cal P}}
\def\R{\overline{R}}
\def\eee{\hfill\rule{1ex}{1ex}}
\def\la{{\lambda}}
\def\La{{\Lambda}}
\def\Lm{{\L_\mu}}
\def\Lla{{\L_\la}}
\def\inv{^{-1}}
\newcommand{\proof}{{\it Proof.\ }}
\newcommand{\bew}{{\it Beweis.\ }}% may  be sc
\newcommand{\va}{\varepsilon}
\def\nml{\normalsize}
\newcommand{\Alg}{\mbox{\rm -Alg}}
\newcommand{\Coalg}{\mbox{\rm -Coalg}}
\newcommand{\Bialg}{\mbox{\rm -Bialg}}
\newcommand{\Hopf}{\mbox{\rm -Hopf}}
\renewcommand{\hom}{\underline{\mbox{\rm hom}}}
\newcommand{\Nat}{\mbox{\rm Nat}}
\newcommand{\rk}{\mbox{\rm rank}}
\newcommand{\tr}{\mbox{\rm tr}}
\renewcommand{\dim}{\mbox{\rm dim}}
\newcommand{\Ext}{\mbox{\rm Ext}}
\newcommand{\Tor}{\mbox{\rm Tor}}

\def\rref#1{(\ref{#1})}
\newcommand{\M}{\mbox{${\cal M}^H$}}
\font\Fraktur=eufm10 scaled\magstep1          % for display- and textstyle
   \newcommand{\fraktur}[1]{\mbox{\Fraktur #1}}  %
   \font\Fraktu=eufm7 scaled\magstep1            % for scriptstyle
   \newcommand{\fraktu}[1]{\mbox{\Fraktu #1}}    %
   \font\Frakt=eufm5 scaled\magstep1             % for scriptscriptstyle
  \newcommand{\frakt}[1]{\mbox{\Frakt #1}}      %
   \def\fr#1{\mathchoice{\fraktur {#1}}            % displaystyle
                        {\fraktur {#1}}            % textstyle
                        {\fraktu {#1}}             % scriptstyle
                        {\frakt {#1}}  }           % scriptscriptstyle

\newcommand{\Ss}{\fr S}

\def\IZ{\mathchoice{ \hbox{${\sf Z}\!\!{\sf Z}$} }
                   { \hbox{${\sf Z}\!\!{\sf Z}$} }
                   { \hbox{$ \scriptstyle {\sf Z}\!\!{\sf Z}$} }
                   { \hbox{$ \scriptscriptstyle {\sf Z}\!\!{\sf Z}$} }  }

\def\db{{\mathchoice{\mbox{\rm db}}
                    {\mbox{\rm db}}
                    {\mbox{\scriptsize\rm db}}
                    {\mbox{\tiny\rm db}} }}

\def\ev{{\mathchoice{\mbox{\rm ev}}
                    {\mbox{\rm ev}}
                    {\mbox{\scriptsize\rm ev}}
                    {\mbox{\tiny\rm ev}} }}

\def\id{{\mathchoice{\mbox{\rm id}}
                    {\mbox{\rm id}}
                    {\mbox{\scriptsize\rm id}}
                    {\mbox{\tiny\rm id}} }}

\def\op{{\mathchoice{\mbox{\rm op}}
                    {\mbox{\rm op}}
                    {\mbox{\scriptsize\rm op}}
                    {\mbox{\tiny\rm op}} }}

\def\ko{{\mathchoice{\mbox{\rm ko}}
                    {\mbox{\rm ko}}
                    {\mbox{\scriptsize\rm ko}}
                    {\mbox{\tiny\rm ko}} }}

\def\ad{{\mathchoice{\mbox{\rm ad}}
                    {\mbox{\rm ad}}
                    {\mbox{\scriptsize\rm ad}}
                    {\mbox{\tiny\rm ad}} }}
\def\sym{{\mathchoice{\mbox{\rm sym}}
                    {\mbox{\rm sym}}
                    {\mbox{\scriptsize\rm sym}}
                    {\mbox{\tiny\rm sym}} }}

\def\R{{\cal R}}
\def\Rbar{\overline{\R}}
\newcommand{\Mod}{\mbox{Mod}}
\newcommand{\Mat}{\mbox{\rm M}}
\newcommand{\Comod}{\mbox{Comod}}

\begin{abstract}{ Bialgebras, defined by means of Yang-Baxter operators which verify the Hecke equation, are considered. It is shown that they are Koszul algebras. Their
Poincar\'e series are calculated via the Poincar\'e series of the corresponding
quantum planes.}\end{abstract}
\newcommand{\Ll}{{\L_\lambda}}

The theory of deformations of function algebras on  algebraic groups was invented by Manin \cite{manin1}, Faddeev, Reshetikhin, Takhtajan \cite{frt}. The idea  is to consider  deformations of coordinate spaces, called quantum planes and then consider the ``symmetries'' of these quantum planes. For example, the quantum group $GL_q(2)$ is the ``symmetry group'' of the quantum planes $A^{2|0}_q$ and $A^{0|2}_q$ \cite{manin1}.

Quantum planes are represented by quadratic algebras, which should be considered as the function algebras on them. Manin considered two quadratic algebras, defined on a finite dimensional vector space $V$ and its dual $V^*$. And in the simplest case, when the commutativity rules are controled by only one element $q$ from the ground field, he obtained the standard deformation of $GL_q(n)$. 

Manin's construction was generalized by Takeuchi \cite{tak90} to the case of orthogonal and symplectic groups. Sudbery \cite{sud} and Mukhin \cite{mukhin} generalized the construction to the case of an arbitrary family of quadratic algebras defined only on $V$. In fact in Manin's construction as well as in the classical case, the R-matrix is symmetric, so we can identify $V$ and $V^*$.

One should mention some earlier works  of Lyubashenko \cite{lyu} and Gurevich \cite{gur}, where the authors  quantized the ``basic symmetry'', i.e. quantized  the category, where we are working in.

Thus being given a family of quadratic algebras, we want to consider their ``symmetries''. In the language of  algebraic geometry, these ``symmetries'' should be represented by a bialgebra (or a Hopf algebra if we need only ``invertible symmetries''), which universally coacts  on the quadratic algebras. It is shown that such a bialgebra (resp. Hopf algebra) always exists.

In this work we consider bialgebras, which are determined as above by only two quadratic algebras, such that the R-matrix obeys the \yb\ equation and Hecke equation.  We call such  bialgebras matrix bialgebras of type $A_n$. We will show that these bialgebras are  Koszul algebras and also consider their Poincar\'e series.

\section{Quadratic algebras and matrix bialgebras}
\subsection{Quadratic algebras} 
Let $k$ be a fixed algebraically  closed field of characteristic zero, $V$ be a $k$-linear vector space of finite dimension $d$. A quadratic algebra (QA)  $A$ over $V$ is defined to be a factor algebra of the tensor algebra over $V$ by an ideal, generated by a set $R(A)$ of quadratic elements, i.e. elements of $V\ot V$. $A$ inherits the grade of $T(V)$ : $A=\oplus_n A_n$. The Poincar\'e series of $A$ is by definition the following formal series:
\bbas P_A(t):=\sum_{n=0}^\infty \dim_kA_nt^n.\eeas
The reader is referred to \cite{manin1} for the definition of the Koszul complex (of second kind) of $A$. $A$ is said to have the Koszul property or to be a Koszul algebra if this Koszul complex is exact \cite{manin1,loef}. The following lemma is due to Backelin \cite{backelin}.
\begin{lem}\label{backelin} A is a Koszul algebra iff the lattice generated by $R(A)$ in $\Vn$ is distributive for all $n\geq 2$.\end{lem}
  A lattice  on $\Vn$ is a set of subspaces of $\Vn$, which is closed under +
  and $\cap$. The lattice generated by $R(A)$   is the one generated by $R_i^n(A), i=1,\cdots ,n-1$, where 
\bbas R_i^n(A):=V^{\ot i-1}\ot <R(A)>_k\ot V^{\ot n-i-1}.\eeas
The lattice $\L$ is distributive if and only if  for all $u,v,w\in \L$  
\bbas u\cap (v+w)=(u\cap v)+(u\cap w),\\
u+(v\cap w)=(u+v)\cap(u+w).\eeas
Note that the two equations are in fact equivalent.
\begin{lem}\label{zerlegung} Let X be a vector space and $\L$ be a lattice on X. Let $X=\bigoplus_{i\in I}X_i$ be a decomposition of X into its subspaces, such  that for all $u\in \L$
\bbas u=\bigoplus_{i\in I}u\cap X_i.\eeas
Then $\L$ is distributive  if and only if for all $i\in I$ 
\bbas \L\cap X_i :=\{u\cap X_i|u\in \L\}\eeas 
 is a distributive lattice on $X_i$.\end{lem}
\proof
Let $u_i:=u\cap X_i$. We have for $u,v\in \L$
\bbas\begin{array}{c}u+v=\oplus(u+v)_i\supset\oplus (u_i+v_i),\\
u\subset \oplus (u_i+v_i), v\subset\oplus (u_i+v_i).\end{array}\eeas Hence $(u+v)_i=u_i+v_i$. Whence the assertion follows.\eee
 \subsection{Matrix bialgebras}\label{section1.2}
 Assume we are given a family of quadratic algebras over a vector space, we want to study
 the algebra $E$, which universally coacts  on this family
 \cite{manin1,ph,sud}.  It is shown that $E$ exists and is a quadratic algebra, and from the universal property it follows that $E$ is a bialgebra. We call such bialgebras  matrix bialgebras.

We consider the following construction. Let $R:V\ot V\lora V\ot V$ be a
diagonalizable operator. Let $R=\sum_{i=1}^{k} c_{i} P_{i}$ be the spectral
decomposition of R. Let $A_i$ be a QA with $R(A_i)=\Im P_i$. Then the matrix bialgebra $E$, determined by the family $\{A_i|i=\dayso{k}\}$, is the factor algebra of the tensor algebra over $E_1:\cong V^*\ot V$ by the ideal, generated by
\bbas R(E)=\Im (\R -1),\eeas
where $\R:=s_{23}(R^{*-1}\ot R)s_{23}:(V^*\ot V)^{\ot 2}\lora (V^*\ot V)^{\ot
  2}$, $s_{23}$ denotes the operator which interchanges elements in 2-nd and
3-rd places of the tensor product. Later on we shall however identify $\R$
with $R^{*-1}\ot R$, acting on $V^{*\ot 2}\ot V^{\ot 2}$. The matrix $R$ is
some times called an R-matrix. The construction is called \yb\ if $R$ obeys the \yb\ equation \cite{mukhin,frt}.

We are interested in the case when $R$ obeys the Hecke  equation
\bbas (R+1)(R-q)=0\eeas
with $q\neq 0$. $R$ is then called  a Hecke operator and  $E$ is  called a
matrix bialgebra of type $A_{d-1}$, i.e. $E$ is considered as a deformation of the
function algebra on the semigroup of $d\times d$ matrices. Thus we have two QA's. We denote them by $\Lambd$ and $S$, where
\bbas R(\Lambd)=\Im (R+1), \ \ R(S)=\Im (R-q).\eeas
 This construction is motivated by the symmetric and antisymmetric tensor algebras over $V$ together with the function algebra on $\End( V)$ coacting upon them.
\subsection{Standard deformations}
 The most interesting example of the construction described in Section  \ref{section1.2} which has many applications in theoretical physics is the case, when $R$ is Drinfeld-Jimbo's R-matrix. We restrict ourselves to the case of R-matrix of series $A_{d-1}$, which gives the standard deformation of $GL(d)$. This deformation plays a crucial role in my work. On some fixed basis $\day{x}{d}$ the operator $R^q$ is given by
\bb R^q:=q\sum_{i=1}^d e_i^i\ot e_i^i+\sqrt{q}\sum_{i,j=1\ i\neq j}e_i^j\ot e_j^i+(q-1)\sum_{i<j}e_i^i\ot e_j^j\label{djmatrix},\ee
where $e_i^j$ is  the operator $e_i^j:V\lora V$, $x_k e_i^j=\delta^j_kx_i$. $R$ obeys the \yb\ and the Hecke equations. The two QA's defined by $R$ are
\bbas & S^q\cong k<x_1,\cdots ,x_d>/(x_ix_j=\sqrt{q}x_jx_i,i< j),\\
\mbox{and}& 
\Lambd^q\cong k<x_1,\cdots ,x_d>/(x_i^2=1,\ x_ix_j=-\sqrt{q}x_jx_i\ ,i\leq j).\eeas
The matrix bialgebra $E$ is defined by:
\begin{eqnarray*}
 z_i^kz_i^l-\sqrt{q}z_i^lz_i^k=0, k<l,\\
z_i^kz_j^k-\sqrt{q}z_j^kz_i^k=0, i<j,\\
z_j^kz_i^l-z_i^lz_j^k=0, i<j,k<l,\\
z_i^kz_j^l+(\sqrt{q}^{-1}-\sqrt{q})z_i^lz_j^k
+z_j^lz_i^k=0,\\ i<j,k<l.
\end{eqnarray*}
It is shown that $S^q, \Lambd^q$ and   $E^q$ have PBW bases \cite{manin1}. By
a theorem of Priddy \cite{priddy}, $S^q ,\Lambd^q$ and   $E^q$ are  then
Koszul algebras.
\subsection{Hecke algebras}\label{hecke}
The reader is referred to \cite{dj1} for a beautiful description of symmetric groups and Hecke algebras. We recall some important facts. The length of an element $w$  of the symmetric group $\Ss_n$ is equal to
$$l(w):=\#\{(i,j)|i< j\ \& \  iw>jw\}.$$
The basic transpositions
$v_{i}=(i,i+1), i=\dayso{n-1}$ generate $\Ss_{n}$. An element $w\in \Ss_{n}$ can be expressed as
a product of $l(w)$ basic transpositions. Let $\B$ denote  the set of the
basic transpositions. 
The elements $v_{i}, i=\dayso{n-2}$ generate a subgroup of
$\Ss_{n}$, isomorphic to $\Ss_{n-1}$. By means of this isomorphism we will
consider  $\Ss_{n-1}$ as a subgroup of $\Ss_{n}$.
The following lemma will be needed in the next section.
\begin{lem}\label{form} For every element $v$ in $\Ss_{n}\setminus\Ss_{n-1}$ there exist elements $w$ and $w'$ in $\Ss_{n-1}$  and $i,j,1\leq i,j\leq n-1$ such that
\bbas & v=v_iv_{i+1}\cdots v_{n-1}w,\\ &  v=w'v_{n-1}v_{n-2}\cdots v_{j}
%\mbox{and}& l(w)=l(w')=l(v)-n+i,
\eeas
for some $1\leq i,j\leq n-1.$\end{lem}
\proof For $w$ in $\Ss_{n-1}$ and $1\leq i\leq n-1$ we have
\bbas  l(v_{i}\cdots v_{n-1}w)= l(w)+n-i\eeas hence for  $w\neq w'$ in $\Ss_{n-1}$  and $1\leq i,j\leq n-1$ 
\bbas
 wv_{i}\cdots v_{n-1}w\neq v_{j}\cdots v_{n-1}w'.\eeas
Whence the first equation of the lemma  follows, the second equation is proved analogously.\eee

The Hecke algebra $\H_n=\H_{q,n}$ as a $k$-space has the basis $T_w,w\in \Ss_n$ with the multiplication subject to:
\begin{enumerate}
\item $T_1=1_{\H_n}$.
\item $T_wT_v=T_{wv}$\ if \ $l(wv)=l(w)+l(v)$.
\item $T_v^2=q+(q-1)T_v$ for $v=(i,i+1), i=\dayso{n-1}.$
\end{enumerate}

We shall allways assume that $[n]_q!\neq 0$ and $q\neq 0$, therefore $\H_n$ is semisimple. If this is the case, $\H_n$ is isormorphic to  a direct product of matrix rings over $k$ \cite{dj1,dj3}. Put 
\bbas x_n=\frac{1}{[n]_q!}\sum_{w\in\Ss_n}T_w,\\
y_n=\frac{1}{[-n]_q!}\sum_{w\in\Ss_n}(-q)^{-l(w)}T_w.\eeas
Then  $x_n,y_n$ are idempotents and 
\bba\label{twxn} T_wx_n=x_nT_w=q^{l(w)}x_n\\
T_wy_n=y_nT_w=(-1)^{l(w)}y_n\eea
(cf. \cite{dj1,dj3}). According to Lemma \ref{form} we have
 \bb\label{xn}[n]_qx_n=x_{n-1}(1+T_{v_{n-1}}+T_{v_{n-1}}T_{v_{n-2}}+\cdots +T_{v_{n-1}}\cdots T_{v_{1}}),\ee
where $v_k$  denotes the basic transposition: $v_k=(k,k+1)$.

 Let us consider the representation $\rho$ of the Hecke algebra $\H_n$ on $\End_k(\Vn)$ induced by a Hecke operator  $R$:
\bbas \rho :\H_n\lora \End_k(\Vn )\label{darstellung},\\
\rho (\Tvi ):x\loma xR^n_i,\ \ i=\dayso{n-1},\eeas
where $R^n_i:=\id^{\ot i-1}\ot R\ot \id^{\ot n-i-1}:\Vn\lora \Vn$.
Let $\Lambd$ and $S$ be quadratic algebras with 
\bbas R(\Lambd)=\Im (R+1),& \ R(S)=\Im (R-q).\eeas
 Then  $\Lambd_n\cong \Im \rho(y_n)\ , \ S_n\cong \Im \rho(x_n)$ \cite{gur}, hence
\bb\label{dimension}s_n:=\dim_kS=\chi (x_n)\ , \ \lambda_n:=\dim_k\Lambd_n= \chi (y_n).\ee
where $\chi$ is the character of $\rho$. Thus the character $\chi$  determines
the Poincar\'e series of $\Lambd$ and $S$. In Section 2.1 I will show that
the Poincar\'e series of $\Lambd$ determines $\chi$.
\section{Matrix bialgebras of type $A_n$}
In this section we define the Schur algebra of a matrix bialgebra of type
$A_{n}$, which will be called R-Schur algebra. We show the double centralizer theorem for it. Using R-Schur algebra one can calculate the Poincar\'e series of the
matrix bialgebra. This will be done in 2.1. In 2.2 we show that the matrix
bialgebra is a Koszul algebra. Using this fact we give in 2.3 the formula for
calculating the Poincar\'e series of the matrix bialgebra via the ones of the
quantum planes.

\subsection{R-Schur algebras and the Poincar\'e series of $E$}
Let $E$ be the matrix bialgebra defined as in the first section by a Hecke
operator $R$. Every homogeneous component $E_{n}$ of $E$ is a finite
dimensional coalgebra, hence its dual $E_{n}^{*}$ is an algebra. The R-Schur
algebra $S_{d}$ on $V$ is defined to be
\bbas S_d=\bigoplus_{n=0}^\infty S_{d,n},& S_{d,n}:=E_n^*,\eeas
where $d$ is the dimension of $V$. Let  $\theta _n $ be the isormorphism
\bbas \theta_n:  \Vsn\ot\Vn \lora E_1^{\ot n}& \\
\theta_{n}(x_1\ot\cdots\ot x_n\ot
y_1\ot\cdots\ot y_n)= & (x_1\ot y_1\ot \cdots\ot x_n\ot y_n) .\eeas 
Let us consider the  natural action of $E_n^*$ on $\Vn$ by 
\bbas \ev_{\Vn}\circ\theta_n^{*}: E^*_n\ot \Vn\lora \Vn\ot\Vsn\ot\Vn\lora \Vn .\eeas  
The Hecke algebra $\H_{n}$ acts on $\Vn$ by the representation $\rho$ in Section \ref{hecke}. 
\begin{edl} The natural action of $E^*_n$ on $\Vn$ and the action $\rho$  of
  $\H_n$ on $\Vn$, induced by $R$, are centralizers of each other in
  $\End_k(V)$.\end{edl}
\proof
Using $\theta_{n}^{*}$ one can identify $E^{*}_{n}$ with a subspace of $\Vn\ot\Vsn$. Then  $E^{*}_{n}$ acts on $\Vn$ as follows. For $f=\sum_i f_i\ot f^*_i \in E^*_n, x\in \Vn$,
$$(\sum_if_i\ot f^*_i)x=\sum_if_i<f^*_i|x>.$$

The Hecke algebra acts on $\Vn$ from the right with
\bbas x\Tvk=xR_k,\ v_k\in \B.\eeas
Since $E^{*}_{n}$ is invariant under $R_{i}$, $k=1,2,\cdots ,n-1$, we have for $f\in E_n^*$
\bbas f=f\R_k^{*-1}=\sum_i( f_i\ot f^*_i)(R_k^{-1}\ot R_k^*)=\sum_i f_iR_k^{-1}\ot f^*_iR^*_k.\eeas
Hence for $k=1,2,\cdots ,n-1$,
 \bbas f(xT_{v_k})=\sum_i f_i<f_i^*|xR_k>=\sum_if_i<f_i^*R_k^*|x>=\\ 
\sum_if_iR_k^{-1}<f_i^*R^*_k|x>R_k=\sum_if_i<f_i^*|x>R_k=(fx)T_{v_k}.\eeas
Thus we obtain the embedding $E_n^*\hookrightarrow \End_{\H_n}(\Vn)$.

Let now $g$  be an element of $ \End_{\H_n}(\Vn)\hookrightarrow \End_k\Vn$.
Using  $\theta$ one can consider  $g$ as an element of $\Vn\ot\Vsn$:
$g=\sum_ig_i\ot g_i^*$. Then $g$ acts on $\Vn$ as follows.
$$gx=\sum_ig_i<g_i^*|x>.$$
Since $g \in \End_{\H_n}(\Vn)$, one has
$$g(x\Tvk)=(gx)\Tvk,\ \forall v_k\in \B ,x\in\Vn .$$
Hence for  $k=1,\cdots ,n-1$,
 $$ \sum_ig_i<g_i^*|x>=\sum_ig_i<g_i^*|x>R_k,\ \forall x\in\Vn ,$$
therefore
 $$ g\R^{*-1}=\sum_i(g_i\ot g_i^*)(R_k^{-1}\ot R_k^*)=\sum_ig_i\ot g_i^*.$$
Thus $E^*_n\cong\End_{\H_n}(\Vn)$ as subalgebras of the algebra $\End_k\Vn$.

Since $\H_{n}$ is semisimple, $E_n^*\cong \End_{\H_n}(\Vn)$ is also semisimple
and by  the density theorem (see for example \cite{cohn} Chapter 10)  the
homomorphism 
\bbas \H_n\lora \End_{E_n^*}\Vn\eeas
is surjective, that completes the proof.\eee

Let $\day{\A}{m}$ be simple subalgebras of $\H_n$. Since the $k$-dimension of
$\A_{i}$ is independent of the field $k$, it is a matrix ring over $k$. Let
$a_{i}$  be the unit element of $\A_{i}$, then $\sum_{i}a_{i}=1_{\H_{n}}$. We have
\bbas  \Vn =\bigoplus_{i=1}^m\Vn a_i=\bigoplus_{i=1}^m\Vn\A_i.\eeas
Hence
\bbas E^{*}_{n}=\End_{\H_n}( \Vn\A_1)\times\cdots\times\End_{\H_n}( \Vn\A_m).\eeas
 $\dim_{k} (\Vn\A_i)$ is equal to $\chi (a_i)$,  where $a_i$ is the unit  element in $\A_i$. Let $r_i$ be the $k$-dimension of the simple $\H_n$-module, which corresponds to $\A_i$. Then
\bbas \End_{\H_n}( \Vn\A_i) =\mbox{M}(\chi ({a_i})/r_i),\eeas  hence
\bb\label{dim} \mbox{dim}_kE_n^*=\sum_i(\frac{\chi (a_i)}{r_i})^2.\ee
Thus, for calculating the dimension of $E_{n}$, it is sufficient to calculate the
character $\chi$ of $\rho$. We show that this can be calculated via the
Poincar\'e series of $S$. We shall do it in three steps.

 Let $C_n:=\{ c_i=v_1v_2\cdots v_{i-1}|v_i=(i,i+1),i=\dayso{n})\}$ be a subset
 of $\Ss_n$. In the first step we show the following lemma.
\begin{lem}\label{lem2.2}For  all $v\in\Ss_n$, $\chi(T_v)$ is a polynomial on
  $\chi(T_{c_i})$ with coefficients being polynomials in q.\end{lem}
\proof 
We use induction. Consider the sequence
 \bbas \Ss_{1}\subset\Ss_{2}\subset\cdots\subset\Ss_{n}\eeas
where $\Ss_{i}$ is generated by $\day{v}{i-1}$ and assume
that the assertion holds for elements of $\Ss_{n-1}$. Let
$v\not\in\Ss_{n-1}$.  According to Lemma \ref{form} one has
$$T_v=T_{\vi} T_{v_{i+1}} \cdots T_{v_{n-1}}T_{v'},\ \  v' \in\Ss_{n-1}.$$ Hence
\bbas\tr R_v=\tr (R_{\vi}\cdots R_{v_{n-1}} R_{v'})=\\
 \tr (R_{v_{n-1}} R_{v'}R_\vi\cdots R_{v_{n-2}}) = \sum_{w\in\Ss_{n-1}}k_w\tr (R_{v_{n-1}} R_w)=\\
 \sum_{w\in\Ss_{n-1}\setminus\Ss_{n-2}}k_w \tr (R_{v_{n-1}} R_w) +\sum_{u\in\Ss_{n-2}}k_u d^{-n}\tr (R_{v_{n-1}} R_u)
\eeas          
where $k_w$ are   polynomials in $q$.
Since $\vn$ and $\Ss_{n-2}$ commute, the terms in the latter sum of the
last part of the above equation are polynomials in  $\tr c_i,i=1\ccc n-1$. Thus it is sufficient to show that $\tr (R_{v_{n-1}} R_w)$ are
polynomials in $\tr c_i, i=1\ccc n-1$ for $w\in\Ss_{n-1}\setminus\Ss_{n-2} $. There exists an element  $w'\in\Ss_{n-2}$ such that $$T_w=T_{v_{k}} T_{v_{k+1}} \cdots T_{v_{n-1}}T_{w'}.$$
Proceed the above process once again, so that we can restrict ourselves to  
showing that $\tr (R_{v_{n-1}} R_{v_{n-2}} R_u), u\in\Ss_{n-2}\setminus\Ss_{n-3}$ are polynomials in  $\tr c_i, i=1\ccc
n$. After $n-2$ times we are led to the element $\tr( R_{v_{n-1}}
R_{v_{n-2}}\cdots  R_{v_1})$. We have $ T_{v_{n-1}}\cdots T_{v_1} =c_n$, which
concludes the proof.\eee

According to the equations (\ref{xn},\ref{twxn}) one has
\bbas [n]_q \chi(x_n)=\tr \rho(x_{n-1} (1+T_{v_{n-1}}+T_{v_{n-1}}T_{v_{n-2}}+\cdots +T_{v_{n-1}}\cdots T_{v_1}))=\\
\tr \rho((1+T_{v_{n-1}}+T_{v_{n-1}}T_{v_{n-2}}+\cdots +T_{v_{n-1}}\cdots T_{v_1}) x_{n-1}  )=\\
\tr \rho (x_{n-1}) + [n-1]_q\tr\rho (T_{v_{n-1}}x_{n-1})=\\
\chi(x_{n-1})+ \tr\rho (T_{v_{n-1}}x_{n-2}(1+T_{v_{n-2}}+\cdots +T_{v_{n-2}}\cdots T_{v_{1}}))=\\
\chi(x_{n-1})+ \tr\rho (T_{v_{n-1}}x_{n-2})+[n-2]_q \tr\rho (T_{v_{n-1}}T_{v_{n-2}}x_{n-1})=\\
\chi(x_{n-1})+d^{-n}\chi(T_{c_2})\chi(x_{n-2})+[n-2]_q \tr\rho (T_{v_{n-1}}T_{v_{n-2}}x_{n-2})=\cdots =\\ \chi(x_{n-1}) +d^{-n}\chi(T_{c_2})\chi(x_{n-2})+
\cdots +d^{-n}\chi(T_{c_{n-1}})\chi(x_{1})+
\chi(T_{c_n}),\eeas
therefore \bb\label{circ}
[n]_q\chi(x_n)=\chi(x_{n-1})+d^{-n}\chi(T_{c_2})\chi(x_{n-2})+\cdots +   d^{-n}\chi (T_{c_{n-1}})\chi(x_{1})+    \chi(T_{c_n}).\ee
That means $\chi (T_{c_i}), c_i\in C_n$ are uniquely determined by $\chi(x_1),\chi(x_2),...,\chi(x_n)$. Consider $c_{k}$ and $x_{k}$ as elements of $\H_{k}$, thus we have
$\chi(x_{k})=s_{k}$. Let denote $p_{0}:=1, p_{k}=\chi_{k}(T_{c_{k+1}}), k\geq
1$ and 
$$P(t)=\sum_{k=0}^\infty p_kt^k.$$
The equation  \rref{circ} is then
\bb\label{p(t)} [n]_qs_n=s_{n-1}p_0+p_1s_{n-2}+\cdots +p_{n-1}s_0\ee
or $P_S(t)\cdot P(t)=\sum_{n=1}^\infty [n]_qs_nt^n.$
 Thus we have shown the following theorem.
\begin{pro}The character $\chi$ of the representation $\rho$ can be determined by rational functions in $\day{s}{n}$ and q.\end{pro}
From this proposition and the equation \rref{dim} the theorem follows immediately.
\begin{edl} Assume that $q\neq 0, [n]_q!\neq 0$ then $\dim_kE_n$ is a rational function on $\day{s}{n}$ and $q$.\end{edl}
%%%%%%%%%%%%%%%%%%%%%%%%%%%%%%%%%%%%%%%%%
%%%%%%%%%%%%%%%%%%%%%%%%%%%%%%%%%%%%%%%%%
\subsection{The Koszul property}
The idea of proving that $S,\Lambd$ and $E$ are Koszul algebras is based on using Lemmas \ref{backelin} and \ref{zerlegung}. We decompose the homogeneous components of those algebras into simple $\H_n$-modules and $\H_{q^{-1},n}\times\H_n$-modules and prove the assertion (in terms of lattices) for each module. The assertion for these  modules is again provided by the Koszul property of the standard deformation and the mentioned above Lemmas.

Let us denote by $\Lambda (d,n)$ the set of all compositions $\lambda=(\day{\lambda}{k}\cdots )$ of $ n$    with $\lambda_a=0$ if $a>d$ and  $\P(n)$ be the set of all partitions of $n$. Then the set of simple $\H_n$-modules can be indexed by $\P(n)$ \cite{dj1}. Let $\L_\mu$, $\mu\in\P(n)$ denote the set of all simple $\H_n$-modules. Then we have for some subset $K$ of $\P(n)$
\bb\label{decomposition} \Vn=\bigoplus_{k\in K}c_k\L_k,\ee
where the integer coefficient $c_k$ denotes that $\L_k$ appears $c_k$ times in the decomposition.
By considering the action of $\H_{q^{-1},n}$ on $\Vsn$ via $R^{*-1}$ we get the similar equation
$$\Vsn =\bigoplus_{j\in J}c'_j\L'_j$$
with $\L'_j$ denoting  simple modules of $\H_{q^{-1},n}$ . Hence 
$$\Vsn\ot\Vn  =\bigoplus_{(j,k)\in J\times K}c_kc'_j\L'_j\ot \L_k$$
is the decomposition of $\Vsn\ot\Vn $ as $\H_{q^{-1},n}\times\H_n$-module. We also have
\bba \Im (R_i+1)=\bigoplus_{k\in K}c_k(\Im(R_i+1)|_{\L_k})\label{lzerl},\\
\Im (R_i-q)=\bigoplus_{k\in K}c_k(\Im(R_i-q)|_{\L_k})\label{szerl},\\
\Im (\R_i-1)=\bigoplus_{(j,k)\in J\times K}c_kc'_j(\Im(\R_i-1)|_{\L'_j\ot \L_k}).\label{ezerl}\eea
We remark that the action of $\H_n$ on $\L_\mu$ does not depend on $R$ but only on $\mu$. 

Let us consider  Drinfeld-Jimbo's matrix $R^q$. From a theorem of Priddy (\cite{priddy}, Theorem 5.2) it follows that  $\Lambd^q, S^q,E^q$ have the Koszul property. To every composition $\lambda$ there corresponds a partition $\overline\lambda $ obtained from $\lambda $ by reordering its elements. The equation (\ref{decomposition}) for $R^q$ is
\bb\label{edpw}\Vn\cong\bigoplus_{\lambda\in\Lambda (d,n)} \L_{\overline\lambda}\ee
  (\cite{dpw}, Proposition 5.1).
If $d>n$, $\P(n)\subset \Lambda (d,n)$, hence all simple modules of $\H_n$
appear in this decomposition. According to Lemma \ref{zerlegung}
$\Im(R_i+1)|_{\L_k},i=\dayso{n-1}$ generate a distributive lattice on
$\L_k$. The same fact is true for $\Im(R_i-q)|_{\L_k},i=\dayso{n-1}$ and
$\Im(\R_i-1)|_{\L'_j\ot\L_k},i=\dayso{n-1}$. Using Lemma \ref{zerlegung} again
we obtain the following theorem.
\begin{edl} Assume that $q\neq 0, [n]_q!\neq 0$ then the algebras $\Lambda,S$ and $E$ are Koszul algebras.\end{edl} 
\subsection{The Poincar\'e series of $E_n$}
\newcommand{\Q}{{\bf Q}}

In this section we will give a formula to calculate the dimension of $E_n$ by the dimensions of $S_i, i=\dayso{n}$. Since $E$ is a Koszul algebra it is enough to calculate the dimension of $B_n$, where $B_n:=\cap_{i=1}^{n-1}\Im (\R_i^n-1)$ \cite{manin1}. We will assume that $q$ is transcendent over $\Q$.

We introduce the following operator $P$: for $S$ being a \yb\ operator on $V\ot V$, $P_n(S)$ operates on $\Vn$:
\bb\begin{array}{l} P_1(S)=\id_V,\nonumber\\
P_n(S)=[n]_q^{-1}(P_{n-1}(S)\ot \id)(\id+S_{n-1}+\cdots +S_{n-1}S_{n-2}\cdots S_1).\end{array}\label{pns}\ee
Let $\Phi_n=P_n(-\R)$ and $\Phi^q_n:=P_n(-\R^q)$, $R^q$ is the Drinfeld-Jimbo's Matrix (\ref{djmatrix}).

For $x=x(q)$ a polynomial on $q$, one considers the following operation $(-)_t$:
\bbas (-)_t:x=x(q)\loma (x)_t:=x(1).\eeas
As we remarked in the previous section, the representation $\rho$ of $\H_n$, induced by $R$ on $\Vn$, restricted on a simple modules $\L_\mu $ does not depend on $R$ any more. The same is true for the representation of $\H_n\ot\H_{q^{-1},n}$ on $E_1^{\ot n}$. For an element $x\in \H_n$ one can then define $(\tr\rho(x)|_\Lm)_t$   and  $(\tr\rho(x))_t:=\sum_{\mu\in I} (\tr\rho(x)|_\Lm)_t$.

Let $P_{n\lambda\mu}:=P_n(-\R)|_{\Ll'\ot \Lm}$, then   $P_{n\lambda\mu}$ does not depend on $R$. 
\begin{edl} Assume that  q is transcendent over  $\Q$ then
\bba \Im\Phi_n=\bigcap_{i=1}^{n-1}\Im(\R^n_i-1)=B_n,\\
 \dim_k \Im\Phi_n=(\tr \Phi_n)_t.\eea
\end{edl}
\proof
The inclusion
\bbas \Im\Phi_n\subset\bigcap\Im(\R^n_i-1)\label{inklusion}\eeas
holds for all Yang-Baxter matrix  $R$. It remains to show the inverse
inclusion. We   
first show it  for the R-matrix  from \rref{djmatrix}.  For  $R=R^q$ with  $q$
transcendent on $\mathbf Q$, $E^q$ has the correct dimension:
\bbas \dim_kE^q_n=\left(\begin{array}{c}d^2+n-1\\ n\end{array}\right) .\eeas
Since $E$ is a Koszul algebra, we have
\bbas  \dim_k\bigcap_i\Im(\R^n_i-1)=\left(\begin{array}{c}d^2\\
    n\end{array}\right).\eeas
Thus
\bb\label{imri}\dim_k\Im\Phi^q_n\leq\left(\begin{array}{c}d^2\\
n\end{array}\right).\ee 
On the other hand, if $q=1$, $\Phi_n^1$ is a projector, hence its eigenvalues are 1 and 0, and one has
\bbas \dim_k\Im\Phi_n^1=\tr\Phi^1_n= \left(\begin{array}{c}d^2\\ n\end{array}\right) .\eeas
Rank of an operator is not less then the number of its non-zero eigenvalues, hence for $q$
transcendent the following holds 
\bbas \dim_k\Im\Phi^q_n\geq (\tr \Phi_n^q)_t=\left(\begin{array}{c}d^2\\ n\end{array}\right) .\eeas
Thus we have
\bb\label{q2.10} \dim_k\Im\Phi^q_n= (\tr \Phi_n^q)_t=\left(\begin{array}{c}d^2\\ n\end{array}\right) \ee
and
\bb\label{q2.11}\Im\Phi^q_n=\bigcap^{n-1}_{i=1}\Im({\R^q}^n_i-1).\ee
The formulas \rref{q2.10} and \rref{q2.11} hold for all $n$ and $q$
transcendent. If $d>n$,  then all simple $\H_n$ modules appear in  the
decomposition \rref{edpw}. Hence 
\bbas \Im\Phi^q_n|_{\Lm}=\bigcap_{i=1}^{n-1}\Im({\R^q}^n_i-1)|_{\Lm},\\
\dim_k\Im \Phi^q_n|_{\Lm}=(\tr\Phi^q_n|_{\Lm})_t.\eeas
Using the Equations \rref{decomposition}--\rref{ezerl} we obtain the assertion
for an arbitrary Hecke operator $R.$\eee

Let $p_k:= (\tr R_{c_{k+1}})_t, k \geq 1$ where $c_k=(\dayso{k})\in \H_k$,
$p_0:=1$. Set 
\bbas P(t)=\sum_{n=0}^{\infty}p_nt^n &\mbox{and}& P_2(t)=\sum_{n=0}^{\infty}p^2_nt^n.\eeas
\begin{edl}Assume that q is transcendent over $\Q$ then
\bb\label{eq212} P_E(t)=exp(\int_0^tP_2(t)),\ee
where $p_n$ can be calculated by the formula
\bb\label{eq213} P(t)=P'_S(t)\cdot P_S(t)^{-1}.\ee\end{edl}
\proof  Let \bbas b_n:=\dim_k\bigcap_{i=1}^{n-1}\Im(\R^n_i-1)=(\tr\Phi_n)_t.\eeas
By definition
\bbas(\tr\R_{c_{k+1}})_t=p_k^2,k\geq 0.\eeas
Since $(\tr(\R_v \R_w))_t=(\tr\R_{vw})_t$ for all $v,w\in\Ss_n$ we have
\bbas n(\tr\Phi_n)_t=(\tr\Phi_{n-1})_t-p_1^2(\tr\Phi_{n-2})_t+\cdots +(-1)^{n-1}p^2_{n-1}\eeas
or
\bbas n b_n=\sum^{n-1}_{k=0}(-1)^kp_k^2b_{n-k-1}.\eeas
Whence the equation  \rref{eq212} follows. On the other hand, from \rref{p(t)}
it follows
\bbas n s_n=\sum^{n-1}_{k=0}p_ks_{n-k-1}.\eeas
Whence the equation \rref{eq213} follows.\eee

 \vspace{1cm}
{\bf Acknowledgement. } The author would like to thank Prof. Yu.I.~Manin for pointing him to these problems and Prof. B.~Pareigis for encouragement and useful discussions.

\end{document}